\documentclass[final]{siamltex}
\usepackage{graphicx} % support the \includegraphics command and options
\usepackage{float}
\usepackage{setspace}
\usepackage{amsmath}
\usepackage{fancyhdr}
\usepackage{multirow}
\usepackage{multicol}
\usepackage{blindtext}
\usepackage{color}
\usepackage{bm}
\usepackage{footmisc}

\title{Localized Patterns in Periodically Forced Systems}

\author{A S Alnahdi\footnote{Department of Applied Mathematics, University of Leeds, Leeds LS2 9JT, UK
(ab.nahdi@gmail.com). Questions, comments, or corrections
to this document may be directed to that email address},  J Niesen, A M Rucklidge and T Wagenknecht\thanks{Deceased.}}

\begin{document}
\bibliographystyle{siam}
\maketitle

\begin{abstract}
Spatially localized, time-periodic structures are common in pattern-forming systems, appearing in fluid mechanics, chemical reactions, and granular media. We examine the existence of oscillatory localized states in a PDE model with single frequency time dependent forcing, introduced in \cite{RS} as phenomenological model of the Faraday wave experiment. In this study, we reduce the PDE model to the forced complex Ginzburg--Landau equation in the limit of weak forcing and weak damping. This allows us to use the known localized solutions found in \cite{BYK}. We reduce the forced complex Ginzburg--Landau equation to the Allen--Cahn equation  near onset, obtaining an asymptotically exact expression for localized solutions. We also extend this analysis to the strong forcing case recovering Allen--Cahn equation directly without the intermediate step. We find excellent agreement between numerical localized solutions of the PDE, localized solutions of the forced complex Ginzburg-Landau equation, and the Allen--Cahn equation. This is the first time that a PDE with time dependent forcing has been reduced to the Allen--Cahn equation, and its localized oscillatory solutions quantitatively studied.\\
This paper is dedicated to the memory of Thomas Wagenknecht.

\end{abstract}


\begin{AMS}\end{AMS}

\pagestyle{myheadings}
\thispagestyle{plain}
\markboth{A S Alnahdi,  J Niesen, A M Rucklidge, and T Wagenknecht}{Localized Patterns in Periodically Forced Systems}

\begin{center}
\noindent\makebox[\linewidth]{\rule{\textwidth}{1pt}} 
%\line(1,0){300}
\end{center}
\section{Introduction}
Localized patterns arise in a wide range of interesting pattern-forming problems. Much progress has been made on steady problems, where bistability between a steady pattern and the zero state leads to localized patterns bounded by stationary fronts between these two states \cite{BK2,BK}. In contrast, oscillons, which are oscillating localized structures in a stationary background, are relatively less well understood \cite{LAF,UMS}. Fluid \cite{AF}, chemical reaction \cite{PQS},  and granular media \cite{UMS} problems have been studied experimentally.  When the surface of the excited system becomes unstable (the Faraday instability),  standing waves are found on the surface of the medium. Oscillons have been found where this primary bifurcation is subcritical \cite{CR}, and these take the form of alternating conical peaks and craters against a  stationary background.

Previous studies have averaged over the fast timescale of the oscillation and have focused on PDE models where the localized solution is effectively steady \cite{AT,BYK,CR}. Here we will seek localized oscillatory states in a PDE with time dependent parametric forcing. We find excellent agreement between oscillons in this PDE and steady structures found in appropriate amplitude equations; this the first complete study of oscillatory localized solutions in a PDE with explicit time dependent forcing.

The complex Ginzburg--Landau (CGL) equation is the normal form description of  pattern forming systems close to a Hopf bifurcation with preferred wavenumber zero  \cite{C}. Adding time dependent forcing to the original problem  results in  a forcing term in the CGL equation, the form of which depends on the ratio between the Hopf and driving frequencies. When the Hopf frequency is half the driving frequency (the usual subharmonic parametric resonance), the resulting PDE is known as the forced complex Ginzburg--Landau  (FCGL) equation:
\begin{equation}\label{eq:byk}
A _{T}=({\tilde{\mu} }+{i\nu})A+(1+i{\kappa})A_{XX}-(1+i{\rho})|A|^{2} A +{\Gamma}\bar{A},
\end{equation}
where all parameters are real, and $\tilde{\mu}$ is the distance from the onset of the oscillatory instability, $\nu$ is the detuning between the Hopf frequency and the driving frequency, $\kappa$ represents the dispersion, $\rho$ is the nonlinear frequency correction, and $\Gamma$ is the forcing amplitude. The complex amplitude, $A(X,T),$ represents the oscillation in a continuous system near a Hopf bifurcation point in one spatial dimension. In the absence of forcing, the state $A=0$ is stable, so $\tilde{\mu}<0$. The amplitude of the response is $|A|,$ and $\arg(A)$ represents the phase difference between the response and the forcing.

\begin{figure}[ht]
\begin{center}
\includegraphics[width=3in]{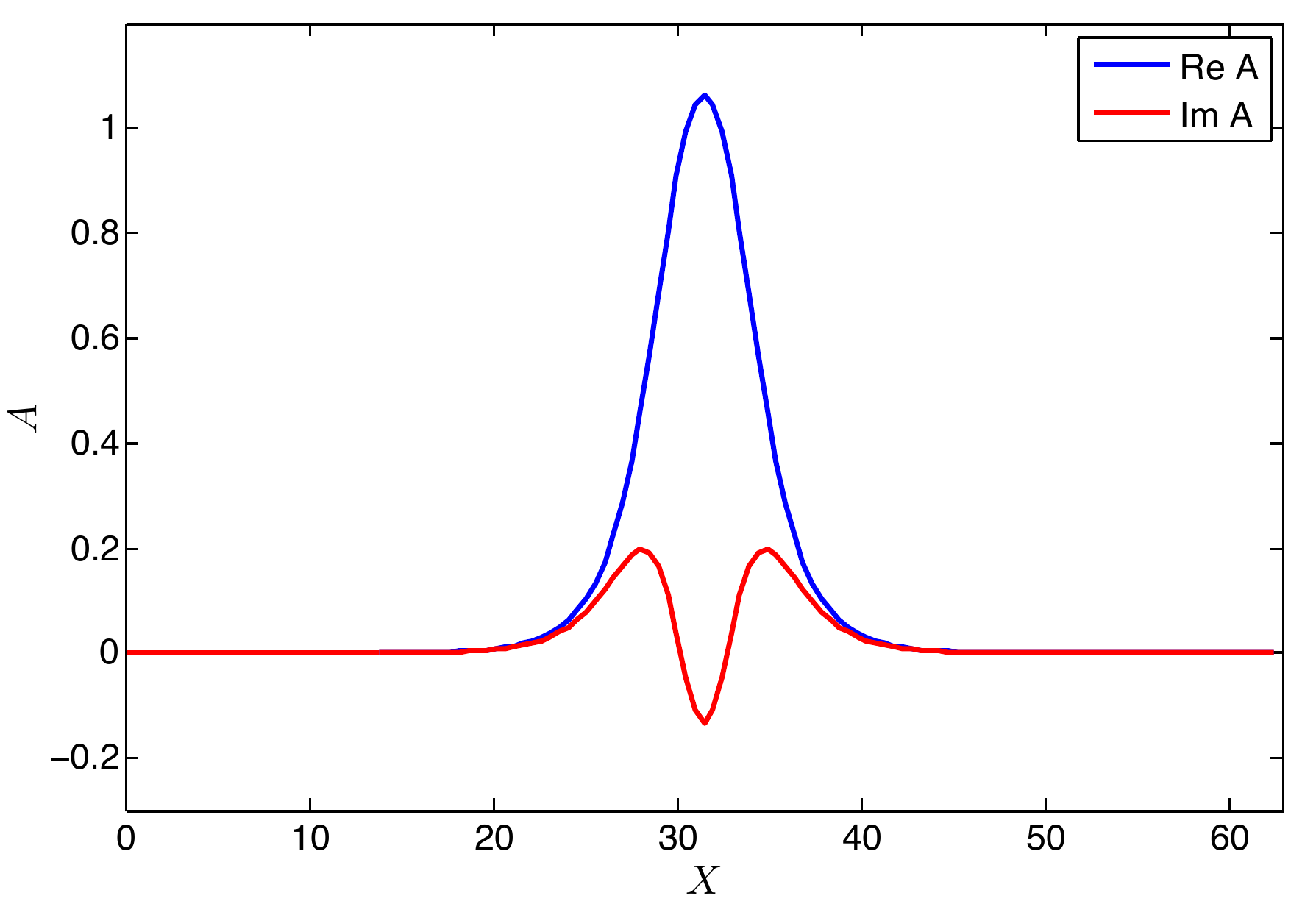}
\end{center}
\caption{Localized solutions of the FCGL equation  (\ref{eq:byk}) with $\tilde{\mu}=-0.5$, $\rho=2.5$, $\nu=2$, $\kappa=-2$, and  $\Gamma=1.496$; the bifurcation point is at $\Gamma_0=2.06$,  following \cite{BYK}. }
\label{fig:fcglp}
\end{figure}

The FCGL equation is a valid description of the full system in the limit of weak forcing, weak damping, small amplitude oscillations and near resonance \cite{CE,EIT}. This model is known to produce localized solutions in 1D \cite{BYK} and in 2D \cite{MS}. It should be noted that these localized solutions have large spatial extent (in the limits mentioned above) and so are different from the oscillons observed in fluid and  granular experiments.
In spite of the cubic coefficient in (\ref{eq:byk}) having negative real part, the initial bifurcation at $\Gamma =\Gamma_0$ is subcritical, the unstable branch turns around in a saddle-node bifurcation, and so there is a nonzero stable solution (the flat state) close to $\Gamma_0$. The localized solution is a homoclinic connection from the zero state back to itself (Figure \ref{fig:fcglp}). Further from $\Gamma_0$, there are fronts (heteroclinic connections) between the zero and the flat state and back.

The aim of this article is to investigate localized solutions in a PDE with parametric forcing, introduced in \cite{RS} as a generic model of parametrically forced systems such as the Faraday wave experiment. We simplify the PDE by removing quadratic terms, by taking the parametric forcing to be $\cos(2t)$, where $t$ is the fast time scale, by working in one rather than two spatial dimensions, and by removing fourth-order spatial derivatives. The resulting model PDE is:
\begin{equation}\label{eq:rs}
U_{t}=({\mu }+{i\omega})U+({\alpha}+i{\beta})U_{xx}
+C|U|^{2}U+i Re(U)F \cos(2t), 
\end{equation}
where $U(x,t)$ is a complex function, $\mu<0$ is the distance from onset of the oscillatory instability, $\omega$, $\alpha$, $\beta$, and $F$ are real parameters, and $C$ is a complex parameter.
\\
In this model the nonlinear terms are chosen to be simple in order that the weakly nonlinear theory and  numerical solutions can be computed easily, and the dispersion relation can be readily controlled. The model shares some important features with the Faraday wave experiment but does not have a clear physical interpretation. The linearized problem reduces to the damped Mathieu equation in the same way that hydrodynamic models of the Faraday instability reduce to this equation in the inviscid limit \cite{BU}. 

We first seek oscillon solutions of (\ref{eq:rs}) by choosing parameter values where (\ref{eq:rs}) can be reduced to the FCGL equation (\ref{eq:byk}). In particular, the preferred wavenumber will be zero, and we will take $F$ to be small, $\mu<0$ to be small, and $\omega$ will be close to 1. We will also consider strong forcing and damping. In the Faraday wave experiment the $k=0$ mode is neutral and cannot be excited, which means experimental oscillons can only be seen with non-zero wavenumbers. This indicates a qualitative difference between this choice of parameters for the PDE model and the Faraday wave experiment.

Here we study equation (\ref{eq:rs}) in two ways. First, in Section 2 we reduce the model PDE asymptotically to an amplitude equation of the form of the FCGL equation (\ref{eq:byk}) by introducing a multiple scales expansion. The numerically computed localized solutions of the FCGL equation (e.g., Figure \ref{fig:fcglp}) will then be a guide to finding localized solutions in the model PDE. Second, we solve the model PDE itself numerically using Fourier spectral methods and Exponential Time Differencing (ETD2) \cite{CM}. We are able to continue the localized solutions using AUTO \cite{AUTO}, and we make quantitative comparisons between localized solutions of the model PDE and the FCGL equation. In Sections 3  and 5 we will do reductions of the FCGL equation and the PDE to the Allen--Cahn equation \cite{AC,YMH} in the weak and strong damping cases respectively; the Allen--Cahn equation has exact localized sech solutions. We give numerical results in Section 4 and conclude in Section 6.

\section{Derivation of the amplitude equation: weak damping case}
In this section we will take the weak forcing, weak damping, weak detuning and small amplitude limit of the model PDE (\ref{eq:rs}), and derive the FCGL equation (\ref{eq:byk}).
Before taking any limits and in the absence of forcing, let us start by linearizing (\ref{eq:rs}) about $U = 0$, and consider solutions of the form 
 $U(x,t)= \exp({\sigma t +ikx})$, where  $\sigma$ is the complex growth rate of a mode with wave number $k$. The growth rate $\sigma$ is given by
\begin{equation}
\sigma=\mu-\alpha k^2 +i (\omega-\beta k^2).
\end{equation}
The forcing $F\cos(2t)$ will drive a subharmonic response with frequency $1$; by choosing $\alpha>0$ and $\omega$ close to 1, we can arrange that a mode with $k$ close to zero will have the largest growth rate. With weak forcing we also need $\mu$, which is negative, to be close to zero, otherwise all modes would be damped. In this case, we are close to the Hopf bifurcation that occurs at $\mu=0$.\\
We now consider the linear theory of the forced model PDE:
\begin{equation}\label{eq:rslin}
U_{t}=({\mu }+{i\omega})U+({\alpha}+i{\beta})U_{xx}
+i Re(U)F \cos(2t), 
\end{equation}
This can be transformed to a Mathieu-like equation \cite{RS}. The normal expectation would be that $\cos(2t)$ would drive a response at frequencies $+1$ and $-1$. However, because $\omega$ is close to $1$, the leading behavior  of (\ref{eq:rslin}) is $$\frac{\partial }{\partial t}U=iU,\qquad\text{or}\qquad \ell_1 U=\left(\frac{\partial }{\partial t}-i\right)U=0.$$
The component of $U$ at frequency $-1$ cancels at leading order, while the component at frequency $+1$ dominates. Furthermore, since $\omega=1+\nu$ with $\nu$ small, and since the strongest response is at or close to wavenumber $k$ where $\omega-\beta k^2=1$, modes with wavenumber $k=0$ will be preferred. Therefore, the leading solution is proportional to $e^{it}$, and so we will seek solutions of the form $U(x,t) =A e^{i t},$ where $A$ is a complex constant. The argument of $A$ relates to the phase difference between the driving force and the response, and is not arbitrary. Later, we will allow $A$ to depend on space and time.\\
To apply standard weakly nonlinear theory, we need the adjoint linear operator $\ell_1 ^\dagger$. First we define an inner product between two functions $f(t)$ and $g(t)$ by
\begin{equation}\label{eq:inprod}
\big\langle f(t),g(t)\big\rangle=\frac{1}{2\pi}\int_0 ^{2\pi} \bar{f}(t) g(t) dt,
\end{equation}
where $\bar{f}$ is the complex conjugate of $f$. With this inner product, the adjoint operator $\ell_1 ^\dagger$, defined by 
$\big\langle f,\ell_1 g\big\rangle =\big\langle \ell_1 ^\dagger f,g\big\rangle$, is given by $$\ell_1 ^\dagger =i-\frac{d }{dt}.$$ The adjoint eigenfunction is then 
$U^\dagger=e^{i t}.$
We take the inner product of (\ref{eq:rslin}) with this adjoint eigenfunction:
\begin{align*}
\begin{split}
0 &=\big\langle U^\dagger, \ell_1 U\big\rangle + \big\langle  U^\dagger, ({\mu }+{i\nu})U+i Re(U)F \cos(2t)\big\rangle \\
&=0+\frac{1}{2\pi}\int_0 ^{2\pi} ({\mu }+{i\nu})U e^{-i t}+\frac{iF}{4} (U+\bar{U})(e^{i t}+e^{-3i t}) dt. 
\end{split}
\end{align*}
We write $U=\sum_{j=-\infty} ^{+\infty} {U_{j} e^{ijt}}$, and $\bar{U}=\sum_{j=-\infty} ^{+\infty}{\bar{U}_{j} e^{-ijt}}$, so
\begin{align*}
\begin{split}
%0&=\frac{1}{2\pi}\int_0 ^{2\pi} ({\mu }+{i\nu})\sum_{j=-\infty} ^{+\infty} {U_{j} e^{i(j-1)t}}\\
%&+\frac{iF}{4} (\sum_{j=-\infty} ^{+\infty} {(U_{j} e^{i(j+1)t}+U_{j} e^{i(j-3)t}+\bar{U}_{j} e^{i(-%j+1)t}+\bar{U}_{j} e^{i(-j-3)t}}) dt\\
0=(\mu+i\nu)U_1 +\frac{iF}{4}(U_{-1}+U_{3}+\bar{U}_1+\bar{U}_{-3}).
\end{split}
\end{align*}
Since the frequency $+1$ component of $U$ dominates at onset, as discussed above, we retain only $U_1$ and $\bar{U}_1$, which satisfy
\[
 \begin{bmatrix}
 \mu+i\nu& \frac{iF}{4} \\
-\frac{iF}{4}& \mu-i\nu\\
 \end{bmatrix}
 \begin{bmatrix}
 U_1\\
  \bar{U}_1
 \end{bmatrix}
=
 \begin{bmatrix}
  0 \\
  0
 \end{bmatrix}\]
This system has a nonzero solution when its determinant is zero; this gives the critical forcing amplitude  $F_0=4{\sqrt{\mu^2+\nu^2}}$. This equation also fixes the phase of $U_1$.\\
To perform the weakly nonlinear calculation, we introduce a small parameter $\epsilon$ and make the substitutions:  $\omega=1+\epsilon^2\nu$,  $F\longrightarrow \epsilon^2 F$,  $\mu\longrightarrow \epsilon^2 \mu$,  and expand the solution $U$ in powers of $\epsilon$ as 
\begin{equation}\label{eq:uexpand}
U=\epsilon U_1+\epsilon^2 U_2 +\epsilon^3 U_3 +...,
\end{equation}
where $U_1$, $U_2$, $U_3$, $....$ are $O(1)$ complex functions.\\
At $O(\epsilon)$, we get $\ell_1 U_1=(\frac{\partial }{\partial t}-i)U_1=0,$ which has solutions of the form
$$U_1=A(X,T) e^{i t},$$ where the amplitude $A$ is $O(1)$, and $X$ and $T$ are slow space and time variables: $T=\epsilon^2 t$, and $X=\epsilon x$. At $O(\epsilon^2),$ we have $U_2(x,t)=0$. At $O(\epsilon^3)$, equation (\ref{eq:rs}) is reduced to
\begin{align*}
\ell_1 U_3+\frac{\partial U_1 }{\partial T}&=( \mu +i \nu)U_1 +(\alpha+ i\beta)\frac{\partial^2 U_1 }{\partial X^2}+C|U_1|^2U_1+i F \cos(2t)Re(U_1),
\end{align*}
We take the inner product with $U^\dagger _1$, and use $\big\langle U ^\dagger _1,\ell_1 U_3\big\rangle=0$ to find the amplitude equation for a long-scale modulation:
\begin{equation}\label{eq:amp1}
A_{T}=( \mu +i \nu)A+ (\alpha +i\beta)A_{XX}+C|A|^2 A +\frac{i F}{4}\bar{A}.
\end{equation} 

We can do a rescaling of the equation (\ref{eq:amp1}) in order to bring it to the  standard FCGL form by rotating $A\longrightarrow A e^{{i\frac{\pi}{4}}}$, which removes the $i$ in front of the $\bar{A}$ term but does not affect any other term. With this, the amplitude equation of the model PDE reads
\begin{equation}\label{eq:amps1}
A_{T}=( \mu +i \nu)A+ (\alpha +i\beta)A_{XX}+C|A|^2 A +\Gamma\bar{A},
\end{equation} 
where $\Gamma=\frac{F}{4}$. A similar calculation in two dimensions yields the same equation but with $A_{XX}$ replaced by $A_{XX}+A_{YY}$.\\
One can see that the amplitude equation (\ref{eq:amps1}) takes the form of the FCGL equation (\ref{eq:byk}). We are now in a position to use the results from \cite{BYK}, where they find localized solutions of  (\ref{eq:byk}), to look for localized solutions of the model PDE (\ref{eq:rs}).

The stationary homogeneous solutions of (\ref{eq:amps1}), which we call the flat states, can easily be computed. These satisfy: 
$$0=( \mu +i\nu)A+C |A|^2 A+\Gamma \bar{A}.$$ 
To solve this steady problem we look for solutions of the form $ A=R e^{i\phi}$, where $R$ is real and $\phi$ is the phase. Dividing by $R e^{i\phi}$ results in:
\begin{equation}\label{eq:steady_sol}
0=( \mu +i\nu)+C R^2  +\Gamma e^{-2i\phi}.
\end{equation} 
We can then separate the real and imaginary parts and eliminate $\phi$ by using  $\sin^2\phi+\cos^2\phi=1$ to get a fourth order polynomial:
\begin{equation}\label{eq:poly}
(C_r^2+C_i^2)R^4+2(\mu C_r+\nu C_i)R^2-\Gamma^2+\mu^2+\nu^2=0,
\end{equation}
where $C=C_r+iC_i$. This can be solved for $R^2$, from which $\phi$ can be determined using (\ref{eq:steady_sol}).\\
Examination of the polynomial ({\ref{eq:poly}) shows that when the forcing amplitude $\Gamma$ reaches $\Gamma_0=\sqrt{\mu^2+\nu^2}$, a subcritical bifurcation occurs provided that $\mu C_r+\nu C_i<0$. A flat state $A^-_{uni} $ is created, which turns into the $A^+_{uni} $ state at $\Gamma_d=\sqrt{-\frac{(\mu C_r+\nu C_i)^2}{(C_r^2+C_i^2)}+\mu^2+\nu^2}$, when a saddle-node bifurcation occurs. We will reduce (\ref{eq:amps1}) further in Section 3 by assuming we are close to onset, and finding explicit expressions for localized solutions.

\section{Reduction to the Allen--Cahn equation: weak damping case}
The FCGL equation (\ref{eq:amps1}) can be reduced to the Allen--Cahn equation \cite[Appendix A]{BYK} by setting $\Gamma=\Gamma_0 +\epsilon_1 ^2\lambda,$ where $\Gamma_0=\sqrt{\mu^2+\nu^2}$ is the critical forcing amplitude, $\lambda$ is the bifurcation parameter, and $\epsilon_1$ is a new small parameter that controls the distance to onset. We expand $A$ in powers of $\epsilon_1$ as
$$A(X,T)=\epsilon_1 A_1(X,T) +\epsilon_1 ^2 A_2(X,T) +\epsilon_1 ^3 A_3 (X,T)+...,$$ where $A_1$, $A_2$, $A_3$ are $O(1)$ complex functions. We further scale $\frac{\partial}{\partial T}$ to be $O(\epsilon_1 ^2)$ and $\frac{\partial}{\partial X}$ to be $O(\epsilon_1$).\\
At $O(\epsilon_1)$ we get
$$0=(\mu+i\nu)A_1+\sqrt{\mu^2+\nu^2} \bar{A}_1, $$
which defines a linear operator \[
\begin{bmatrix}
\mu+i\nu& \sqrt{\mu^2+\nu^2}\\
\sqrt{\mu^2+\nu^2}& \mu-i\nu\\
\end{bmatrix}
\begin{bmatrix}
A_1\\
\bar{A}_1
\end{bmatrix}
=
\begin{bmatrix}
0 \\
0
\end{bmatrix}\]
The solution is $A_1=B(X,T) e^{i\phi_1}$, where $B$ is real, and the phase $\phi_1$ is fixed by $ e^{-2i\phi_1}=-\frac{\mu+i\nu}{\sqrt{\mu^2+\nu^2}}$. This gives
$$\phi_1=\tan^{-1} \left(\frac{\mu+\sqrt{\mu^2-\nu^2}}{\nu}\right),$$ 
At $O(\epsilon_1 ^3)$, we have 
\begin{equation}\label{eq:real_amp}
B_T e^{i\phi_1}= (\mu+i\nu) A_3 +(\alpha+i\beta) B_{XX} e^{i\phi_1}+C B^3 e^{i\phi_1}+\lambda Be^{-i\phi_1}+\Gamma_0 \bar{A}_3.
\end{equation}
We take the complex conjugate of (\ref{eq:real_amp}) and multiply this by $e^{-i\phi_1}$, and then add  (\ref{eq:real_amp}) multiplied by $e^{i\phi_1}$ to eliminate $A_3$. With this, equation (\ref{eq:amps1}) reduces to the Allen--Cahn  equation 
\begin{equation}\label{eq:real_amp2}
B_T= \frac{-\lambda \sqrt{\mu^2+\nu^2}}{\mu}B+\frac{(\alpha \mu+\beta \nu)}{\mu} B_{XX}+\frac{\mu C_r+\nu C_i}{\mu} B^3.
\end{equation}
We can readily find localized solutions of (\ref{eq:real_amp2}) in terms of hyperbolic functions. This leads to an approximate oscillon solution of (\ref{eq:amps1}) of the form 
\begin{equation}\label{eq:exact_sol1}
A=\sqrt{\frac{2(\Gamma-\Gamma_0)\sqrt{\mu^2+\nu^2}}{\mu C_r+\nu C_i}}\text{sech}\left(\sqrt{\frac{(\Gamma-\Gamma_0)\sqrt{\mu^2+\nu^2}}{(\alpha\mu+\beta\nu)}}X\right)e^{i\phi_1},
\end{equation}
provided  $\Gamma< \Gamma_0$, $\mu <0$, $\mu C_r+\nu C_i < 0$, and $\alpha \mu+ \beta \nu <0$.
Note that in the PDE (\ref{eq:rs}) we have the assumption $U_1=\epsilon A e^{it}$, therefore the spatially localized oscillon is given approximately by
\begin{equation}\label{eq:exact_sol2}
U_{loc}=\sqrt{\frac{(F-F_0)\sqrt{\mu^2+\nu^2}}{2(\mu C_r+\nu C_i)}}\text{sech}\left(\sqrt{\frac{(F-F_0)\sqrt{\mu^2+\nu^2}}{4(\alpha\mu+\beta\nu)}} x \right)e^{i(t+\phi_1)},
\end{equation}
again provided $F <F_0$.
We compare the approximate solution $U_{loc}$ with a numerical solution of the PDE below, as a dotted line in Figure \ref{fig:loc_curve_auto1}(a).

\section{Numerical results: weak damping case}

In this section, we present numerical solutions of (\ref{eq:byk}) (in the form written in (\ref{eq:amps1})) and (\ref{eq:rs}), using the known \cite{BYK} localized solutions of (\ref{eq:byk}) to help find similar solutions of (\ref{eq:rs}), and comparing the bifurcation diagrams of the two cases.

We use both time-stepping methods and continuation on both PDEs. For time-stepping, we use a pseudospectral method, using FFTs with  up to 1280 Fourier modes, and  the exponential time differencing method ETD2 \cite{CM}, which has the advantage of solving the non-time dependent linear parts of the PDEs exactly. We treat the forcing term ($\Gamma \bar{A}$ and $ Re(U) \cos(2t)$) with the nonlinear terms.

For continuation, we use AUTO \cite{AUTO}, treating $x$ as the time-like independent variable, to find steady solutions of the FCGL (\ref{eq:amps1}). For the PDE (\ref{eq:rs}), we represent solutions with a truncated Fourier series in time with the frequencies $-3$, $-1$, $1$ and $3$. The choice of these frequencies comes from the forcing $Re(e^{it})\cos(2t)$ in the PDE, taking  $U=e^{it}$ as the basic  solution, as described above. \\
Following \cite{BYK} we will take illustrative parameter values for the amplitude equation (\ref{eq:amps1}): $\mu=-0.5,$ $ \alpha=1$, $\beta=-2$, and $C=-1-2.5i$, and solve the equation on domains of size $L_X= 20\pi$. For (\ref{eq:rs}), we use $\epsilon=0.1$, which implies $\mu=-0.005$, $F=0.04\Gamma$,  $\omega=1.02$, $L_x=200\pi$, and use the same $\alpha$, $\beta$, and $C$. We show examples of localized solutions in the FCGL equation and the PDE (\ref{eq:rs}) in Figure \ref{fig:examples_pde_fcgl}, demonstrating the quantitative agreement as expected between the two.
\begin{figure}[h!]
  \centering
  \begin{tabular}{  c  c  }
    \begin{minipage}{.45\textwidth}
      \includegraphics[trim=4cm 7cm 2cm 7cm, clip=true, totalheight=0.3\textheight]{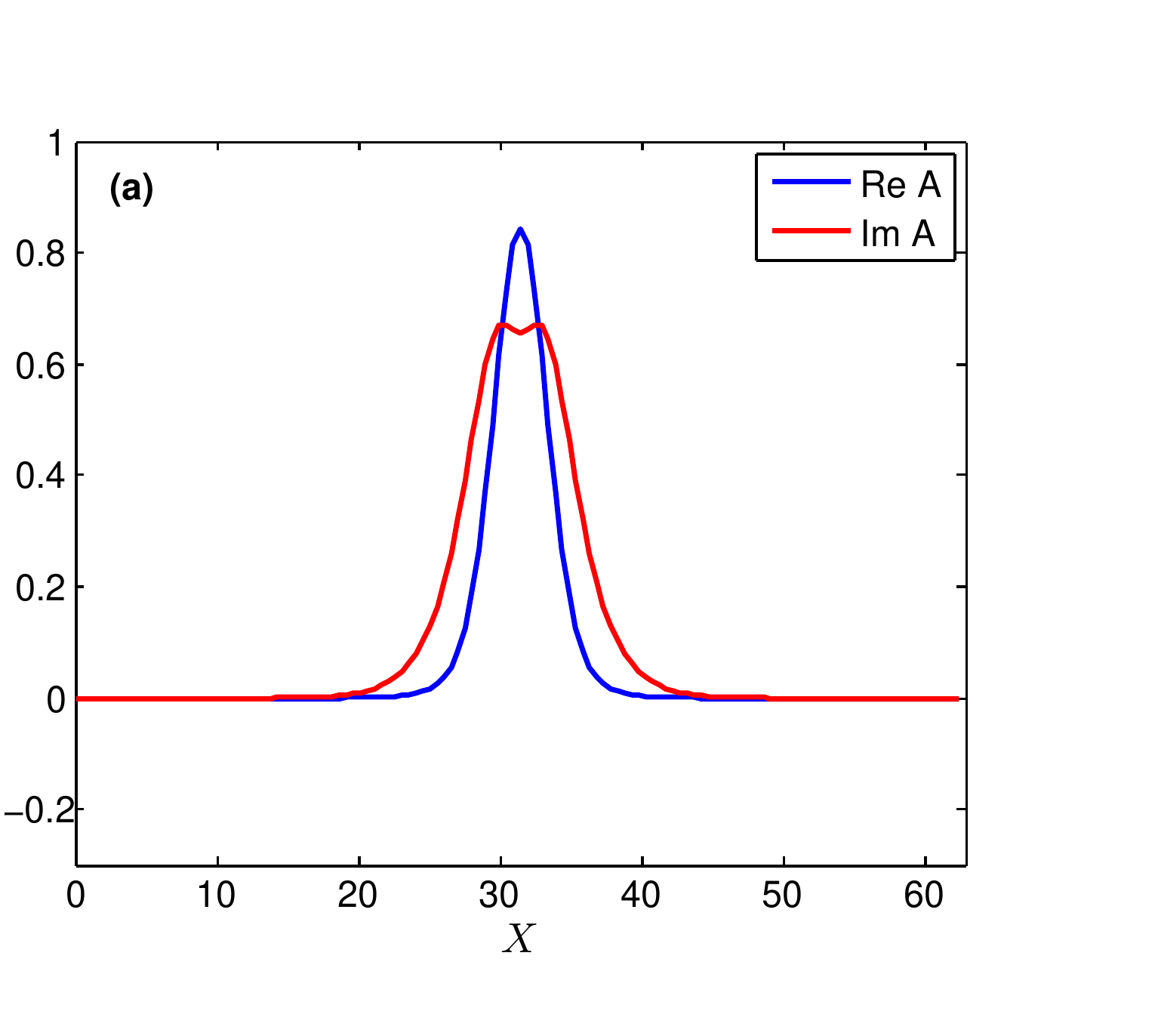}
    \end{minipage}
    &
    \begin{minipage}{.50\textwidth}
\includegraphics[trim=3cm 7cm 2cm 7cm, clip=true, totalheight=0.3\textheight]{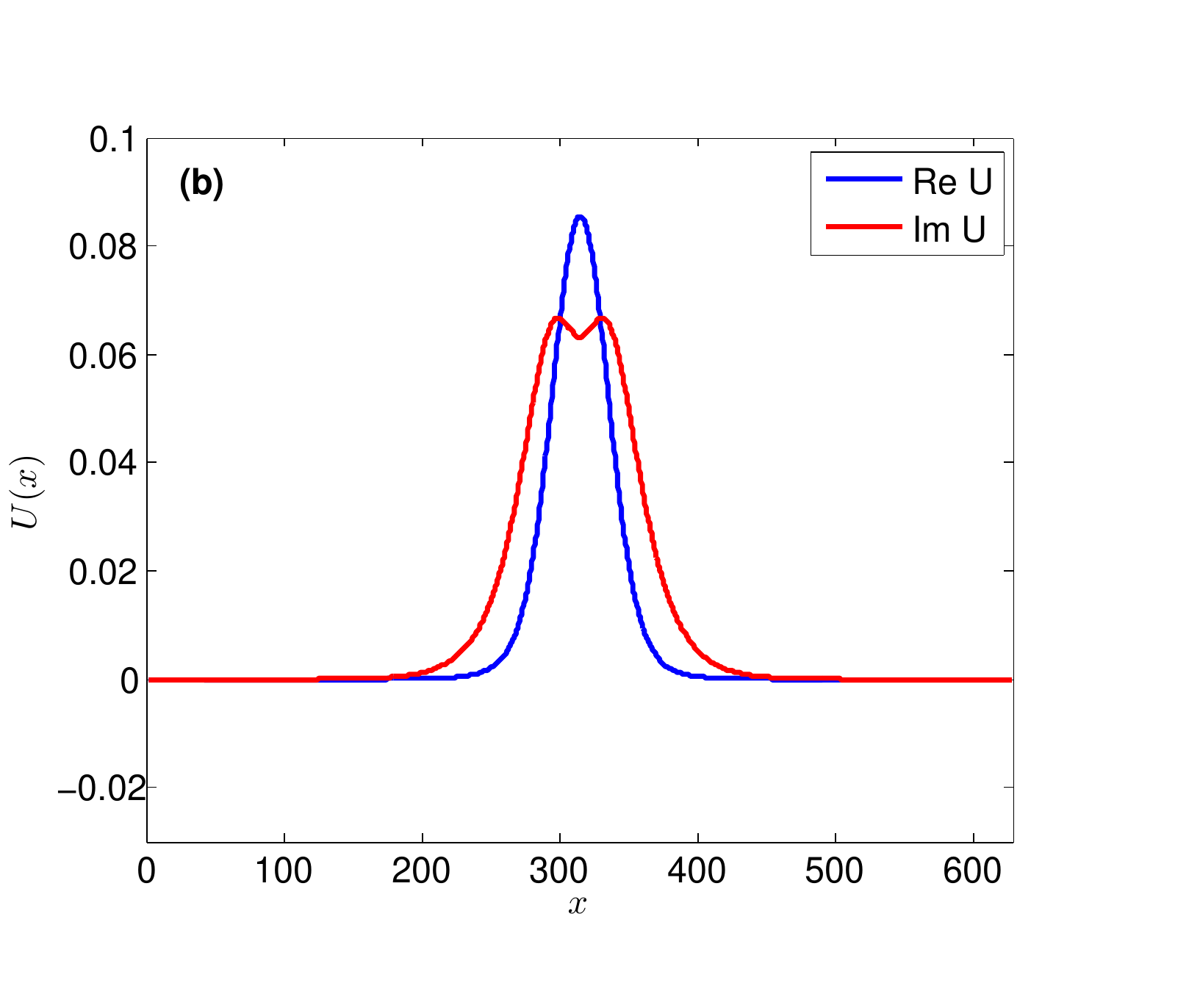} 
\end{minipage}
  \end{tabular}
\caption{(a) Example of a localized solution to the FCGL equation (\ref{eq:amp1}) with $\mu=-0.5$, and $F=5.984$. (b) Example of a localized solution to the PDE model (\ref{eq:rs}) with $\mu=-0.5\epsilon^2$, and $F=5.984\epsilon^2$, where $\epsilon=0.1$. In both models $\alpha=1$, $\beta=-2$,  and  $\nu =2$, and $C=-1-2.5i$. Note the factor of $\epsilon$ in the scalings of the two axes. }
\label{fig:examples_pde_fcgl}
\end{figure}\\
In all bifurcation diagrams we present solutions in terms of their norms \\
$$N=\sqrt{\frac{2}{L_x}\int^{L_x}_0 |U|^2\,dx},$$
We computed (following \cite{BYK}) the location of these stable localized solutions  in the ($\nu$,$\Gamma$) parameter plane, shown in green  in Figure \ref{fig:nu_ga_loc}.  In this figure one can see that the region of localized solutions starts where $\mu C_r+\nu C_i=0$, when the primary bifurcation changes from supercritical to subcritical \cite{DMC,KC}, and gets wider as $\nu$ increases. We also show the bistability region of the amplitude equation between the primary $(\Gamma_0)$ and the saddle-node $(\Gamma_d)$ bifurcations.\\
\begin{figure}[h!]
\begin{center}
\includegraphics[width=3.5in, height=2.5in]{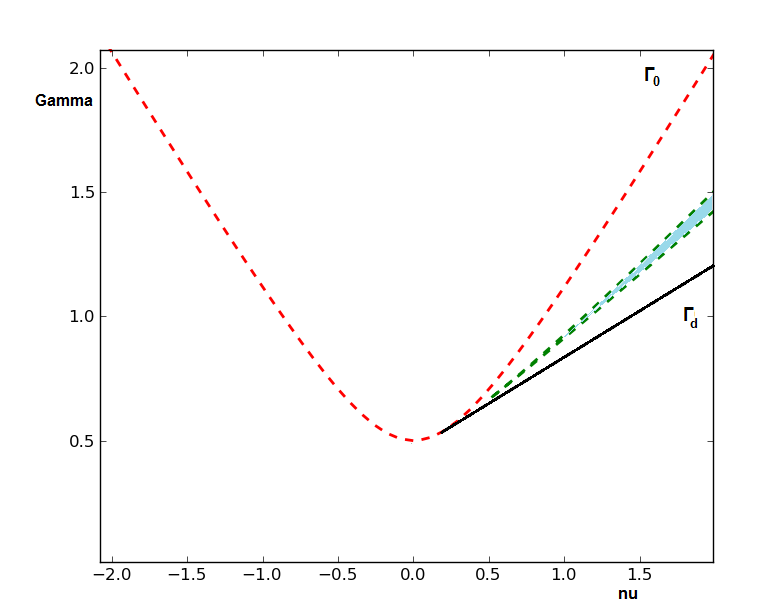}
\end{center}
\caption{The $(\nu,\Gamma)$-parameter plane for FCGL equation (\ref{eq:amps1}), $\mu=-0.5$,  $\alpha=1$, $\beta=-2$, and $C=-1-2.5i$, recomputed following \cite{BYK}. Stable localized solutions exist in the shaded green region. The dashed red line is the primary pitchfork bifurcation at $\Gamma_0=\sqrt{\mu^2+\nu^2}$, and the solid black line is the saddle-node bifurcation at $\Gamma_d$.}
\label{fig:nu_ga_loc}
\end{figure}

\begin{figure}[h!]
\begin{center}
\includegraphics[width=3.5in, height=3in]{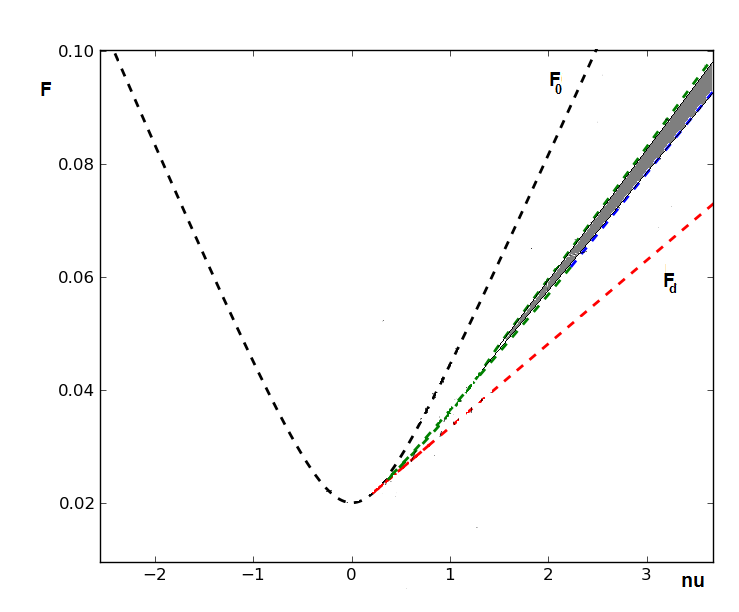}
\end{center}
\caption{The $(\nu,F)$-parameter plane of the PDE model (\ref{eq:rs}) with $\mu=-0.005$, $\alpha=1$, $\beta=-2$, and $C=-1-2.5i$. Stable localized solutions exist in the shaded grey region. The dashed black line is the primary pitchfork bifurcation and the dashed red line is the saddle-node bifurcation at $F_d$.}
\label{fig:F_nu_curve_pde}
\end{figure}

Part of the difficulty of computing localized solutions in the PDE comes from finding parameter values where these are stable.  In the FCGL equation with $\nu=2$, stable localized solutions occur between $\Gamma_{1}^*=1.4272$ and $\Gamma_{2}^*=1.5069$. In the PDE with parameter values as above, we therefore estimate that the stable localized solutions should exist between $F^{*}_{1}=0.04\Gamma_{1} ^{*}=0.0573$ and  $F^{*}_{2}=0.0600$. We found by time-stepping a stable oscillatory spatially localized solution in the PDE model (\ref{eq:rs}) at $F=0.058$ and used this as a starting point for continuation with AUTO. We found stable localized solutions between saddle-node bifurcations, at $F^*_1=0.05688$ and  $F^*_2=0.06001$, which compares well with the prediction from the FCGL equation. In addition, the bistability region was determined by  time-stepping to be between $F_{d}=0.04817$ and $F_{0}=0.08165$. As $\nu$ is varied, the  grey shaded region in Figure \ref{fig:F_nu_curve_pde} shows the region where stable localized solutions exist in the PDE.

The snaking regions of the PDE model and the FCGL equation are presented in Figure \ref{fig:loc_curve_rs4mod}. In this figure we rescale the PDE, so we can plot the bifurcation diagrams of the amplitude equation and the PDE model in top of each other. The agreement is excellent. Examples of localized solutions are given in Figure \ref{fig:loc_curve_auto1} (a)-(f) as we go along the localization curve.  Our comparison between  results from the FCGL equation (\ref{eq:amps1}) in Figure  \ref{fig:nu_ga_loc} and results from  the model PDE (\ref{eq:rs}) in Figure \ref{fig:F_nu_curve_pde} shows excellent agreement. 

Note the decaying spatial oscillations close to the flat state in Figure \ref{fig:loc_curve_auto1} (c)-(f): it is these that provide the pinning necessary to have parameter intervals of localized solutions. These parameter intervals become narrower as the localized flat state becomes wider (see Figure \ref{fig:loc_curve_rs4mod}) since the oscillations decay in space, in contrast with the localized solutions found in the subcritical Swift--Hohenberg equation \cite{BK2}.

In this study so far our calculations have been based on assuming weak damping and weak forcing. Next, we study the PDE in the strong forcing case. %{ref any one explain this in details}.
\begin{figure}[h!]
\begin{center}

\includegraphics[trim=3cm 7cm 2cm 7cm, clip=true, totalheight=0.5\textheight]{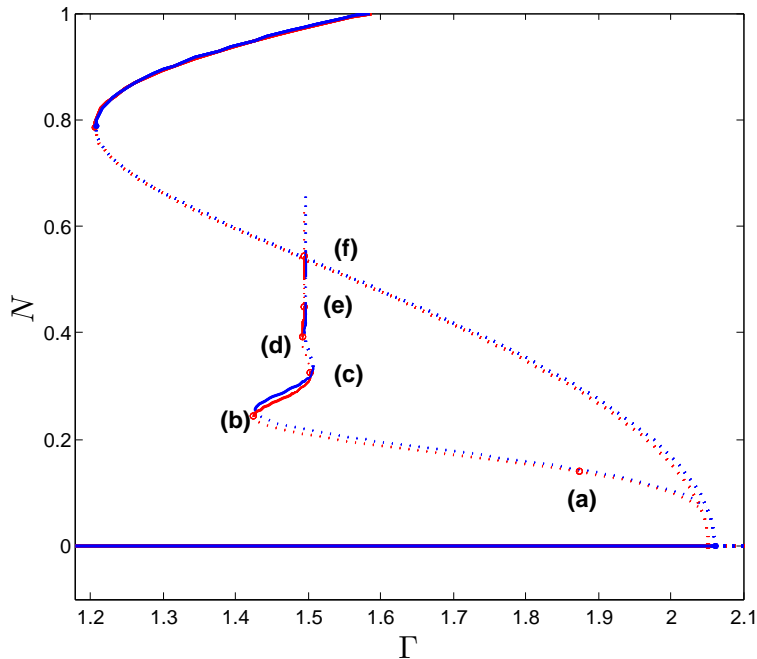}
\end{center}
\caption{The red curves correspond to bifurcation diagram of the PDE model and the blue curves  correspond to the FCGL equation. Solid (dashed) lines correspond to stable (unstable) solutions. For the PDE we use $F=4\epsilon^2 \Gamma$. Parameters are otherwise as in Figure \ref{fig:examples_pde_fcgl}. Example solutions at the points labeled (a)-(f) are in Figure \ref{fig:loc_curve_auto1}. Bifurcation point in the FCGL is $\Gamma_{0}=2.06$, and in the PDE is   $\Gamma_{0}=2.05$.}
\label{fig:loc_curve_rs4mod}
\end{figure}

\begin{figure}[h!]

  \centering
  \begin{tabular}{  c  c  }
    
    \begin{minipage}{.45\textwidth}
    
\includegraphics[trim=1cm 7cm 2cm 7cm, clip=true, totalheight=0.25\textheight]{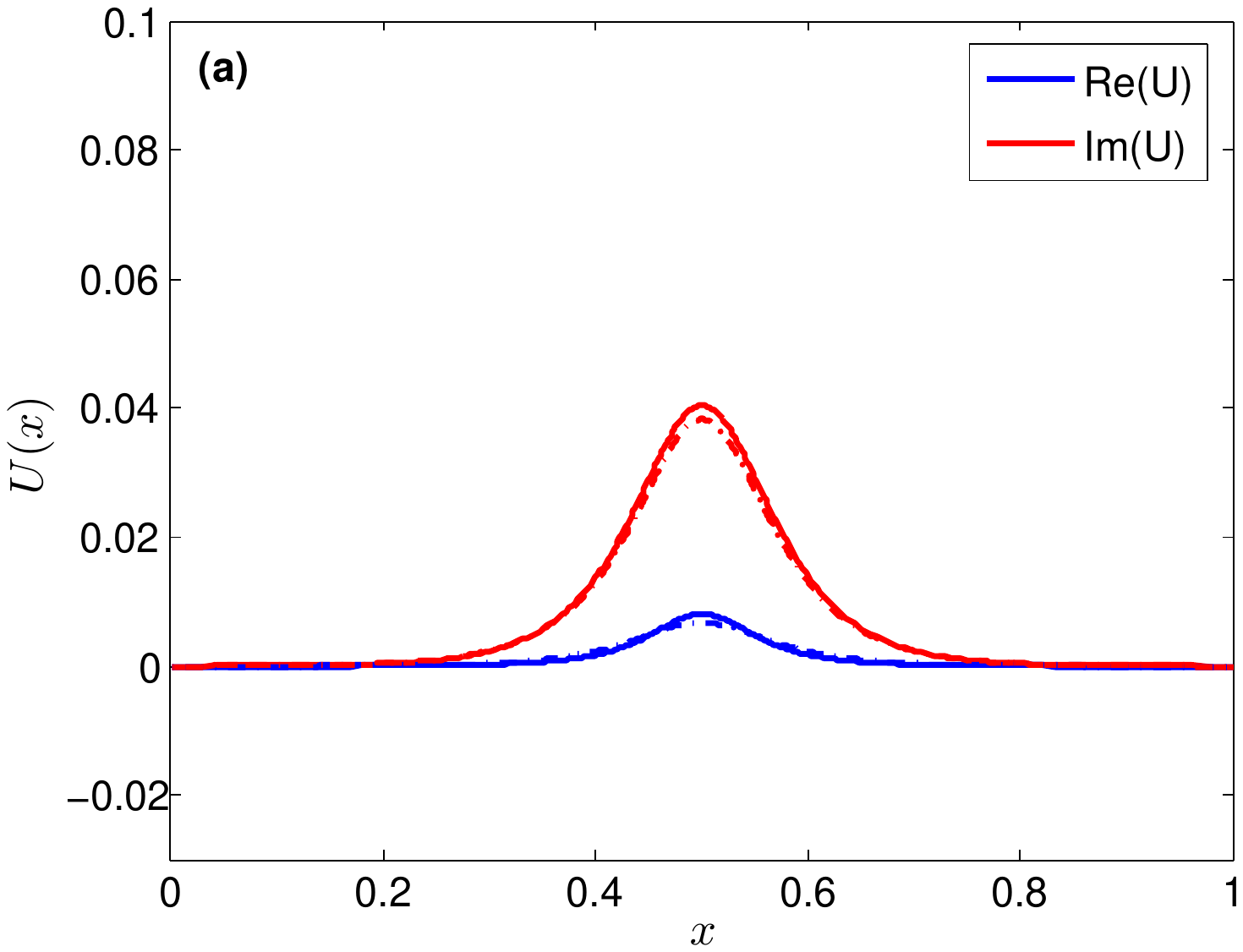}\\
      \includegraphics[trim=1cm 7cm 2cm 7cm, clip=true, totalheight=0.25\textheight]{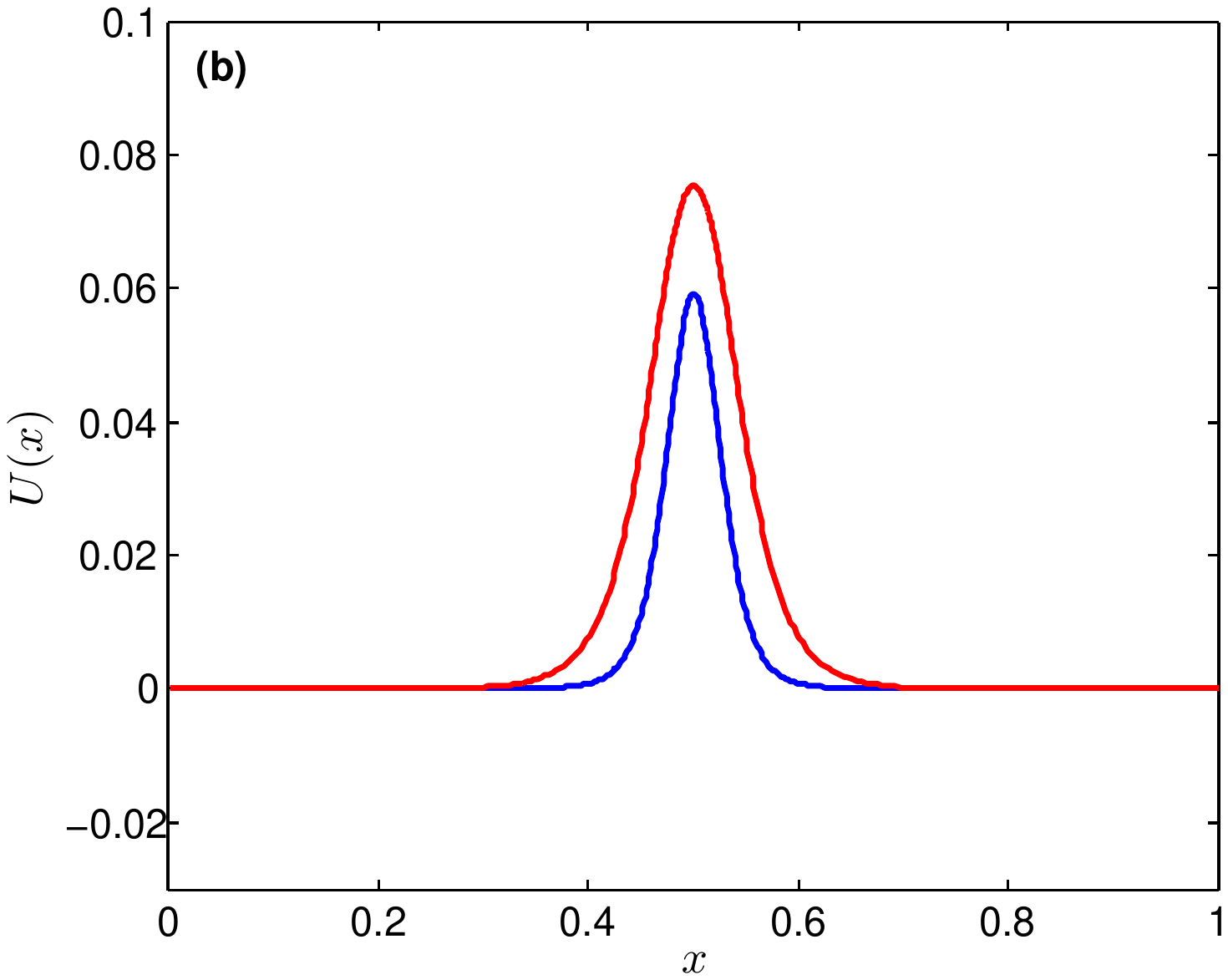}\\
         \includegraphics[trim=1cm 7cm 2cm 7cm, clip=true, totalheight=0.25\textheight]{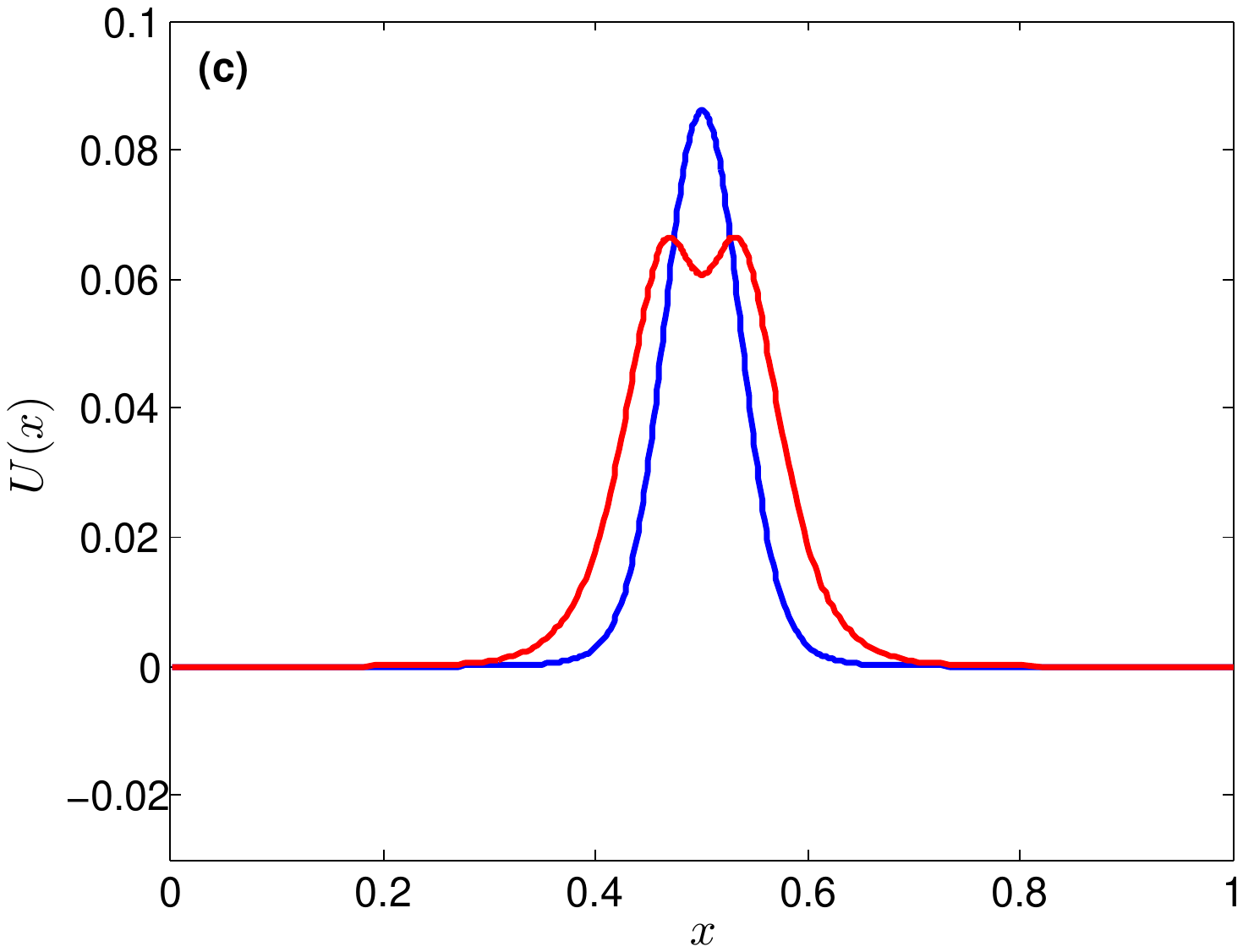}

    \end{minipage}
    &
    \begin{minipage}{.4\textwidth}
           \includegraphics[trim=1cm 7cm 2cm 7cm, clip=true, totalheight=0.25\textheight]{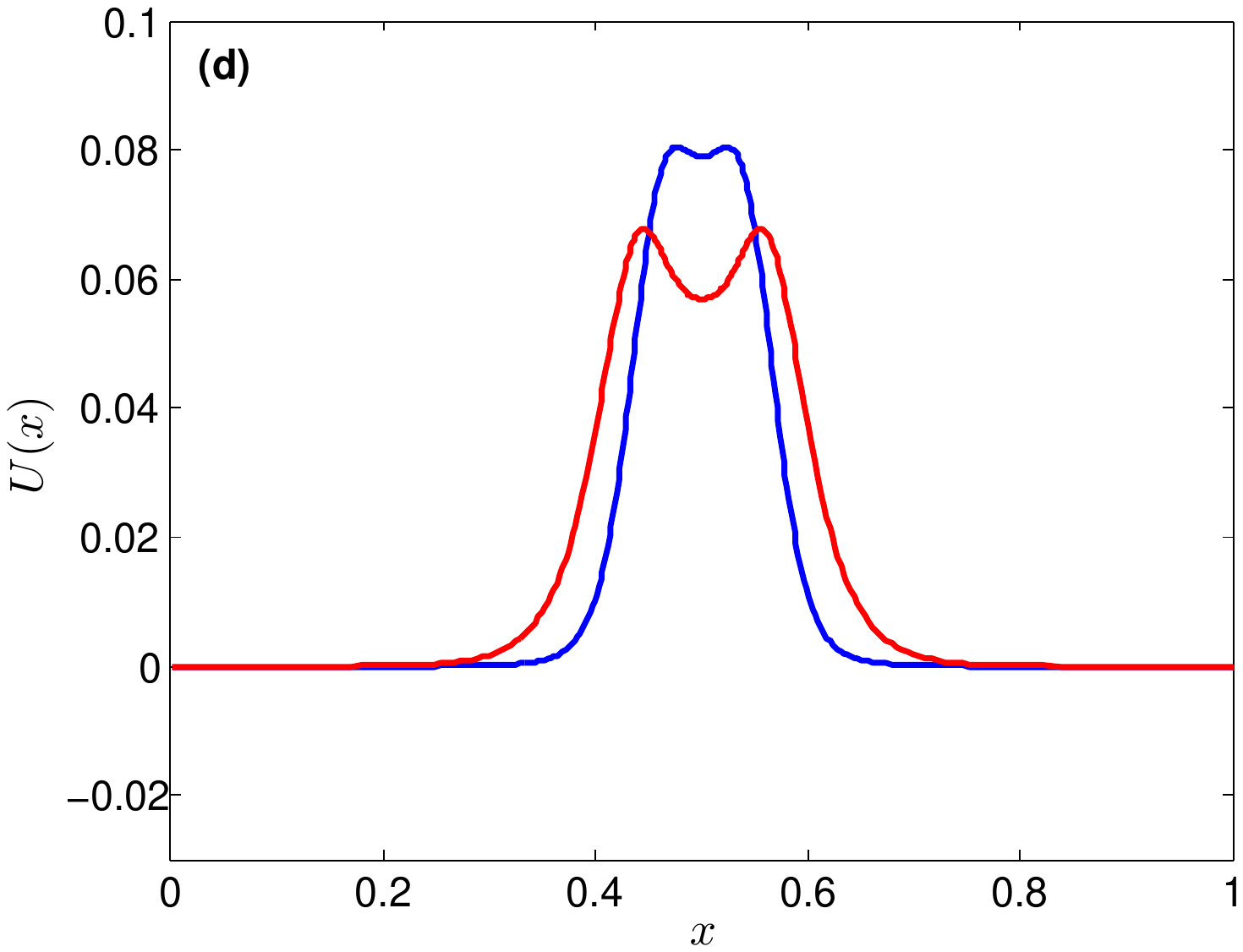}\\
          \includegraphics[trim=1cm 7cm 2cm 7cm, clip=true, totalheight=0.25\textheight]{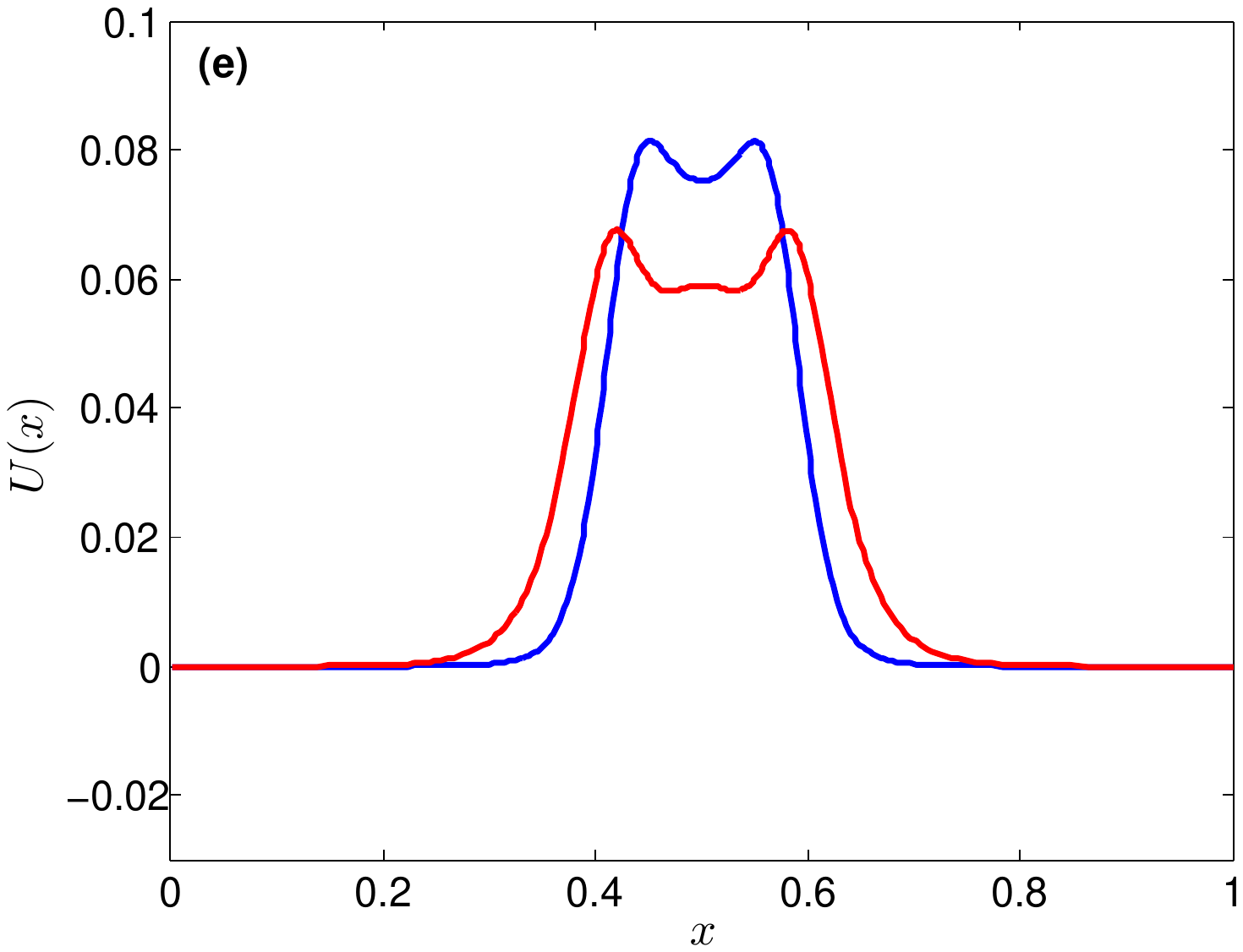} \\
              \includegraphics[trim=1cm 7cm 2cm 7cm, clip=true, totalheight=0.25\textheight]{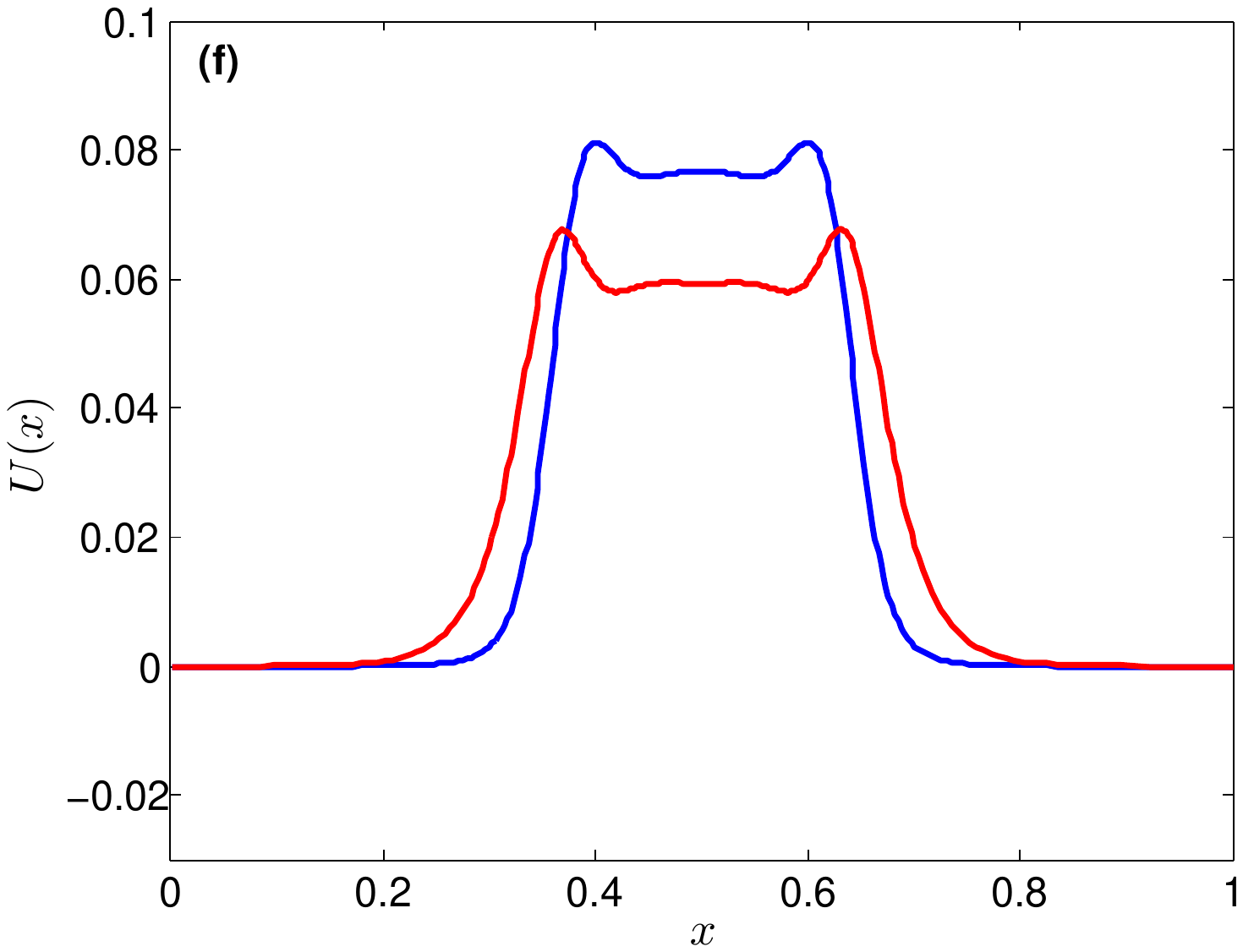} 

\end{minipage}

  \end{tabular}

\caption{Examples of solutions to  (\ref{eq:rs}) equation along the localized branch with $\mu=-0.005$, $\alpha=1$, $\beta=-2$, $\nu=2$, and $C=-1-2.5i$. Bistability region is between $F_0=0.08165$ and $F_d=0.048173$, and localized oscillons branch is between ${F_1}^*=0.05688$ and ${F_2}^*=0.06001$. (a) $F=0.07499$.  (b) $F=0.05699$. (c)$F=0.06015$. (d) $F=0.05961$. (e) $F=0.05976$. (f) $F=0.05975$. Dot lines represent the real (blue) and imaginary (red) parts of $U_{loc}$.}

\label{fig:loc_curve_auto1}
\end{figure}

\section{Reduction of the PDE to the Allen--Cahn equation: strong damping case}
In the strong damping, strong forcing case, the linear part of the PDE is not solved approximately by $U_1= e^{i t}$. Rather, a Mathieu equation must be solved numerically to get the eigenfunction \cite{RS}. In this case, weakly nonlinear calculations lead to the Allen--Cahn equation directly, without the intermediate step of  the FCGL equation (\ref{eq:byk}) with its $\Gamma \bar{A}$ forcing. The advantages of reducing the PDE to the Allen--Cahn equation are that localized solutions in this equation are known analytically, and that demonstrates directly the existence of localized solutions in the PDE model.\\
We write the solution as  $U=u+iv$,
where $u(x,t)$ and $v(x,t)$ are real functions. Thus, equation ({\ref{eq:rs}) is written in terms of real and imaginary parts of $U$ as 
\begin{equation}\label{eq:reu}
\begin{aligned}
\frac{\partial u}{\partial t}&=(\mu+\alpha \nabla^{2})u-(\omega+\beta \nabla^{2})v+C_r(u^2+v^2)u-C_i(u^2+v^2)v,\\
\frac{\partial v}{\partial t}&=(\omega+\beta \nabla^{2})u+(\mu+\alpha \nabla^{2})v+C_r(u^2+v^2)v+C_i(u^2+v^2)u+f(t)u.
\end{aligned}
\end{equation}
We begin our analysis by linearizing (\ref{eq:reu})  about $u=0$ and $v=0$. We write the  periodic forcing function as $ f(t)=f_c(t)(1+\epsilon_1^{2}\lambda)$, where $f_c(t)=F_c\cos(2t)$. Here, $F_c$ is the critical forcing amplitude, which must be determined numerically, and is where the trivial solution loses stability. We seek a critical eigenfunction of the form 
\begin{equation}\label{eq:solution}
U=p_1(t)+iq_1(t),
\end{equation}
where $p_1(t)$ and $q_1(t)$ are  real $2\pi$-periodic functions. Note that in writing $u+iv$ in this form, we are taking the critical wavenumber to be zero. The analysis follows that presented in \cite{RS}, but in the current work the spatial scaling and the chosen solution are different, again because the critical wavenumber is zero. Substituting into (\ref{eq:reu}) at onset leads to
\begin{equation}\label{eq:systempq}
\begin{aligned}
\left[\frac{\partial}{\partial t}-\mu \right]p_1&=-\omega q_1, \\
\left[\frac{\partial}{\partial t}-\mu \right]q_1&=\omega p_1+f_c(t)p_1,
\end{aligned}
 \end{equation}
which can be combined to give a damped Mathieu equation 
\begin{align*}
\left[\frac{d}{dt}-\mu\right]^2 p_1+\left(\omega ^2+f_c(t)\omega\right) p_1=0,
\end{align*}
or
\begin{equation}\label{eq:math1}
\ddot{p}_1 -2\mu\dot{p}_1+\left(\mu^2+\omega^2+f_c(t)\omega\right)p_1=L p=0,
\end{equation}
defining a linear operator 
$$L=\frac{\partial^2 }{\partial t^2}-2\mu\frac{\partial }{\partial t}+(\mu^2+\omega^2+\omega f_c(t))$$
The critical forcing function $f_c(t)=F_c\cos(2t)$ is determined by the condition that (\ref{eq:math1}) should have a non-zero solution $p_1(t)$, from which $q_1(t)$ is found by solving the top line in (\ref{eq:systempq}).
Using the inner product (\ref{eq:inprod}), we have the adjoint linear operator, given by 
$$L^\dagger=\frac{\partial ^2}{\partial t^2}+2\mu\frac{\partial }{\partial t}+(\mu^2+\omega^2+\omega f_c(t)).$$
The adjoint equation is  $L^\dagger p^\dagger _1=0$,  where $p ^\dagger _1$ is the adjoint eigenfunction, which is computed numerically.

In order to reduce the model PDE (\ref{eq:rs}) to the Allen--Cahn equation, we expand solutions in powers of $\epsilon_1$ as \\
\begin{equation}\label{eq:expu}
\begin{aligned}
u=\epsilon_1 u_1+\epsilon_1^2 u_2 +\epsilon_1^3 u_3 +...,\\
v=\epsilon_1 v_1+\epsilon_1^2 v_2 +\epsilon_1^3 v_3 +...,
\end{aligned}
\end{equation}
where $\epsilon_1\ll 1$ and $u_1$, $u_2$, $u_3$, $...$, $v_1$, $v_2$, $v_3$, $...$ are $O(1)$ real functions. We introduce the slow time variable $T=\epsilon_1 ^2t$ and the slow space variable $X=\epsilon_1 x$. Substituting  equation (\ref{eq:expu})  into (\ref{eq:reu}), the associated equations at each power of $\epsilon_1$ are as follows. 
At $O(\epsilon_1)$,  the linear argument above arises, and we have $u_1+i v_1=B(X,T)(p_1+i q_1)$, where $p_1+iq_1$ is the  critical eigenfunction, normalized so that $\langle p_1+i q_1, p_1+i q_1\rangle=1$, and $B$ is a real function of $X$ and $T$. Note that the phase of the response is determined by the critical eigenfunction.
At $O(\epsilon_1 ^3)$, the problem is written as
\begin{align*}
\left(\frac{\partial }{\partial t}-\mu\right) u_3+\frac{\partial u_1}{\partial T}&=-\omega v_3+\alpha \frac{\partial^2 u_1}{\partial X^2}-\beta\frac{\partial^2 v_1}{\partial X^2}+C_{r}(u_1^2+v_1^2)u_1-C_i(u_1^2+v_1^2)v_1,\\
\left(\frac{\partial }{\partial t}-\mu\right)v_3+\frac{\partial v_1}{\partial T}&=\omega u_3+f_c(t)u_3+\lambda f_c(t)u_1+\alpha \frac{\partial^2 v_1}{\partial X^2}+\beta\frac{\partial^2 u_1}{\partial X^2}+C_{r}(u_1^2+v_1^2)v_1\\
&\quad{}+C_i(u_1^2+v_1^2)u_1.
\end{align*}
Eliminating $v_3$, we find
\begin{equation}\label{eq:real_amp_pde}
\begin{aligned}
L u_3&=-\left(\frac{\partial }{\partial t}-\mu\right) \frac{\partial u_1}{\partial T}+\omega\frac{\partial v_1}{\partial T} \\
&\quad{}+\left(\frac{\partial }{\partial t}-\mu\right) \left(\alpha \frac{\partial^2 u_1}{\partial X^2}-\beta\frac{\partial^2 v_1}{\partial X^2}\right)\\
&\quad{}-\omega\left(\alpha \frac{\partial^2 v_1}{\partial X^2}+\beta\frac{\partial^2 u_1}{\partial X^2}\right)-\omega\lambda f_c(t)u_1\\
&\quad{}-\omega\left( C_{r}\left(u_1^2+v_1^2\right)v_1+C_i\left(u_1^2+v_1^2\right)u_1\right)\\
&\quad{}+\left(\frac{\partial }{\partial t}-\mu\right) \left(C_{r}\left(u_1^2+v_1^2\right)u_1-C_i\left(u_1^2+v_1^2\right)v_1\right).
\end{aligned}
\end{equation}
We apply the solvability condition to equation (\ref{eq:real_amp_pde}) $ \langle p ^\dagger _1, L u_3\rangle=0$. We substitute the solution $u_1=B p_1$, and $v_1=B q_1$ into equation (\ref{eq:real_amp_pde}), and then we take the inner product between ${p}^\dagger _1$ and  this equation. Note that we use $\left(\frac{\partial }{\partial t}-\mu\right) p_1=-\omega q_1$, so the equation can be then written as
\begin{equation}
\begin{aligned}
\left\langle p ^\dagger _1,2\left(\frac{\partial }{\partial t}-\mu\right)  p_1\right \rangle\frac{\partial B}{\partial T}&=-\left\langle p ^\dagger _1,\omega f_c(t)p_1\right\rangle \lambda B\\
&\quad{}+\left\langle p^\dagger _1, \left(\left(\frac{\partial }{\partial t}-\mu\right)  (\alpha p_1-\beta q_1)-\omega\left(\alpha q_1+\beta p_1\right)\right)\frac{\partial^2 B}{\partial X^2}\right\rangle\\
&\quad{}+\bigg\langle p^\dagger _1, -\omega\left( C_{r} \left(p_1^2+q_1^2\right)q_1+C_i\left(p_1^2+q_1^2\right)p_1\right)\\
&\qquad{ }+\left(\frac{\partial }{\partial t}-\mu\right) \left(C_{r}\left(p_1^2+q_1^2 \right)p_1-C_i \left(p_1^2+q_1^2 \right)q_1\right)\bigg\rangle B^3,
\end{aligned}
\end{equation}
We find coefficients of the above equation by computing the inner products numerically. Therefore, the PDE is reduced to the Allen--Cahn equation as
\begin{equation}\label{eq:amp}
B_T=1.5687 \lambda B+ 11.1591 B_{XX}+ 9.4717B^3,
\end{equation}
for the parameter values in Figure \ref{fig:strong_example} (a).
Note that $U=\epsilon_1 U_1$,  $X=\epsilon_1 x$, and $\epsilon_1^2 \lambda=\frac{F}{F_0}-1,$
so that the spatially localized solution takes the form
\begin{equation}
 U_{loc}=\sqrt\frac{{-3.1374(\frac{F}{F_0}-1)}}{9.4717} \text{sech}\left( \sqrt\frac{{-1.5687(\frac{F}{F_0}-1)}}{11.1591}  x \right) \left( p_1(t)+i q_1(t) \right).
\end{equation} 
Thus, we have found approximate examples of localized solutions of the PDE, which are qualitatively similar to those found in the weak damping case. Figure \ref{fig:strong_example} (a) shows the comparison between the numerical solution and $U_{loc}$. This solution is continued using AUTO to compute a bifurcation diagram in Figure \ref{fig:strong_example}, again qualitatively similar to the weak damping case. 
\begin{figure}[H]
  \centering
  \begin{tabular}{  c  c  }
    \begin{minipage}{.4\textwidth}
      \includegraphics[trim=.01cm 9cm 3cm 9cm, clip=true, totalheight=0.23\textheight]{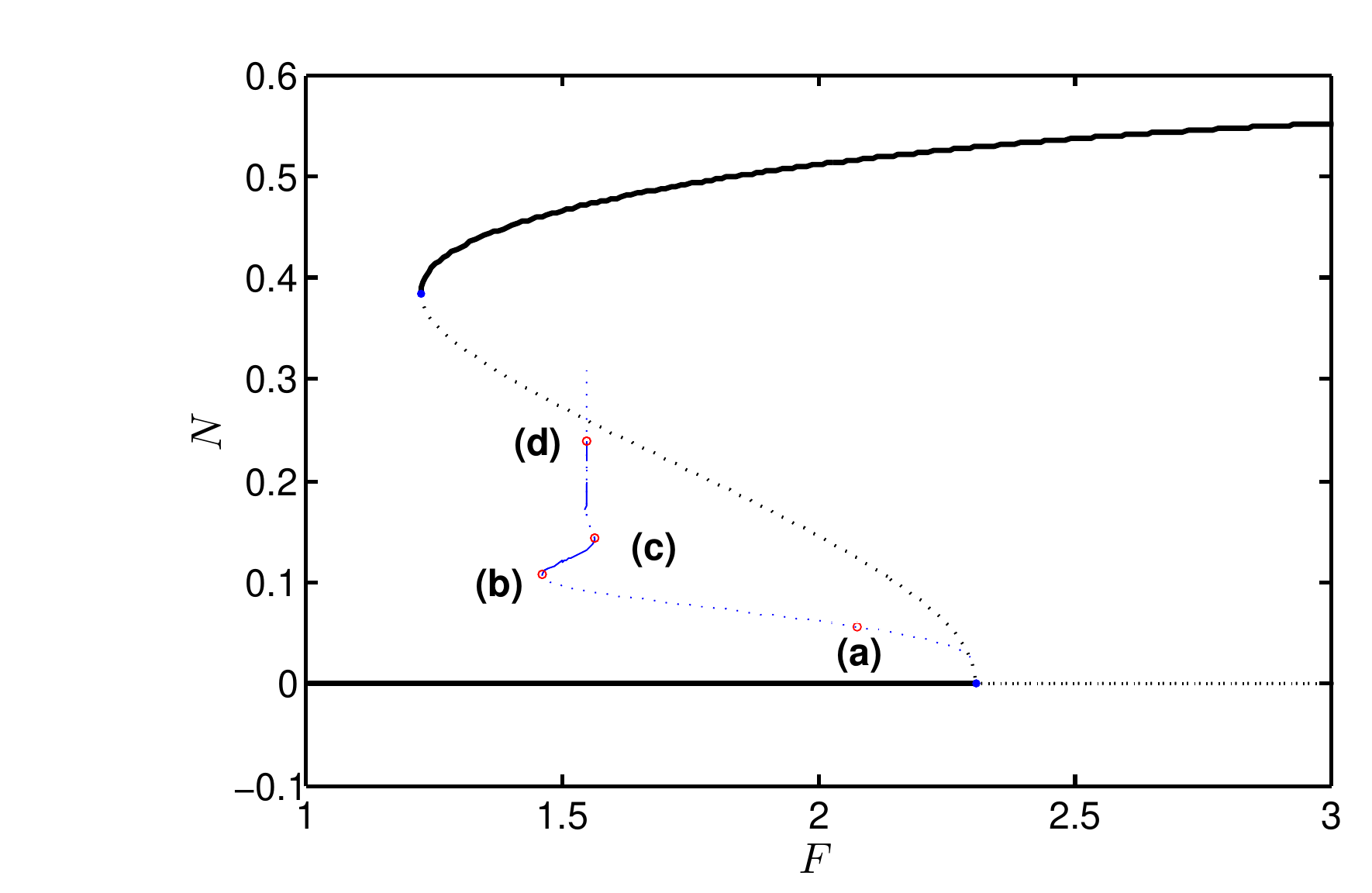}\\
    \end{minipage}
    \\
    \begin{minipage}{.4\textwidth}
           \includegraphics[trim=.2cm 9cm 2cm 7cm, clip=true, totalheight=0.2\textheight]{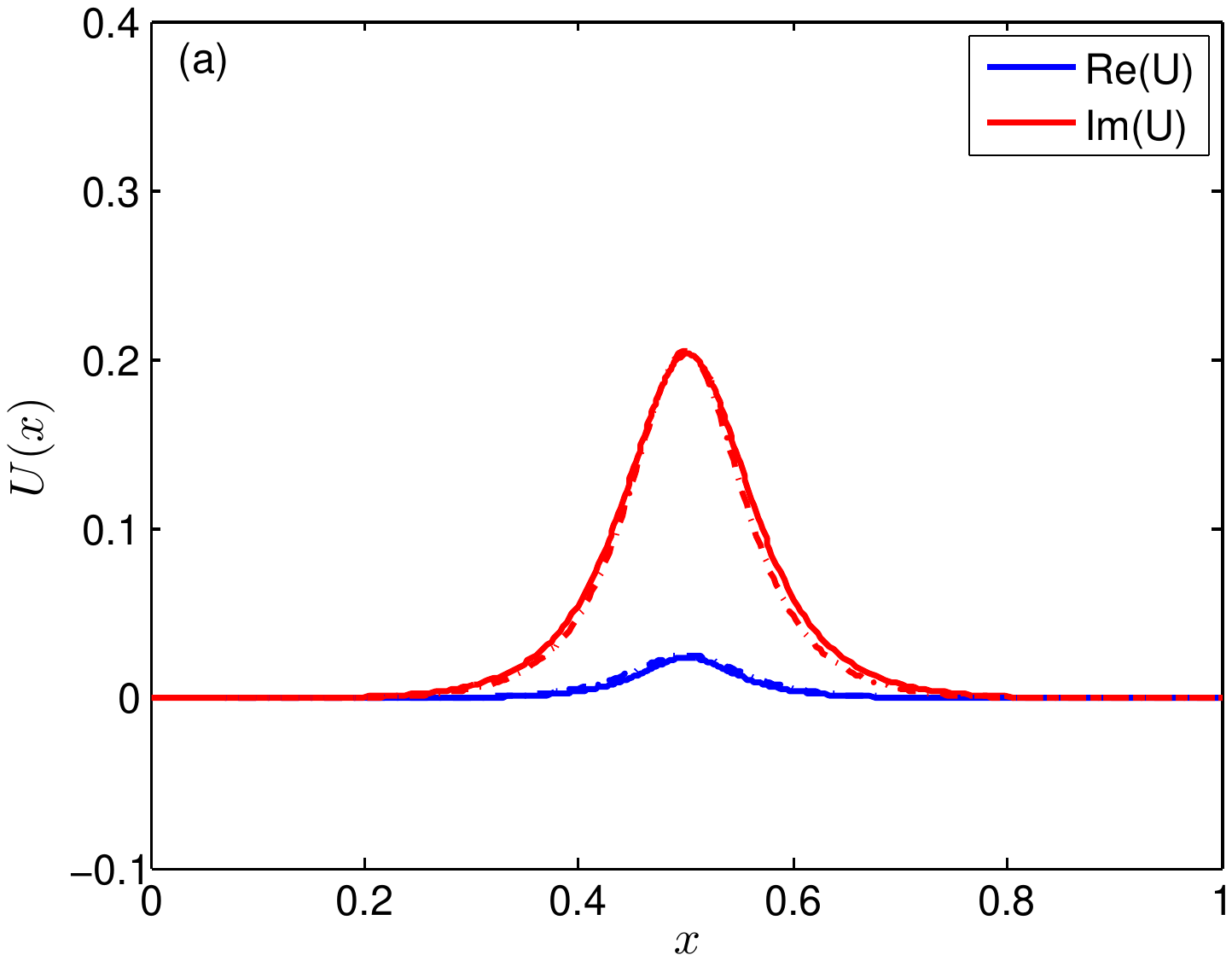}\\
 \includegraphics[trim=.2cm 9cm 2cm 7cm, clip=true, totalheight=0.2\textheight]{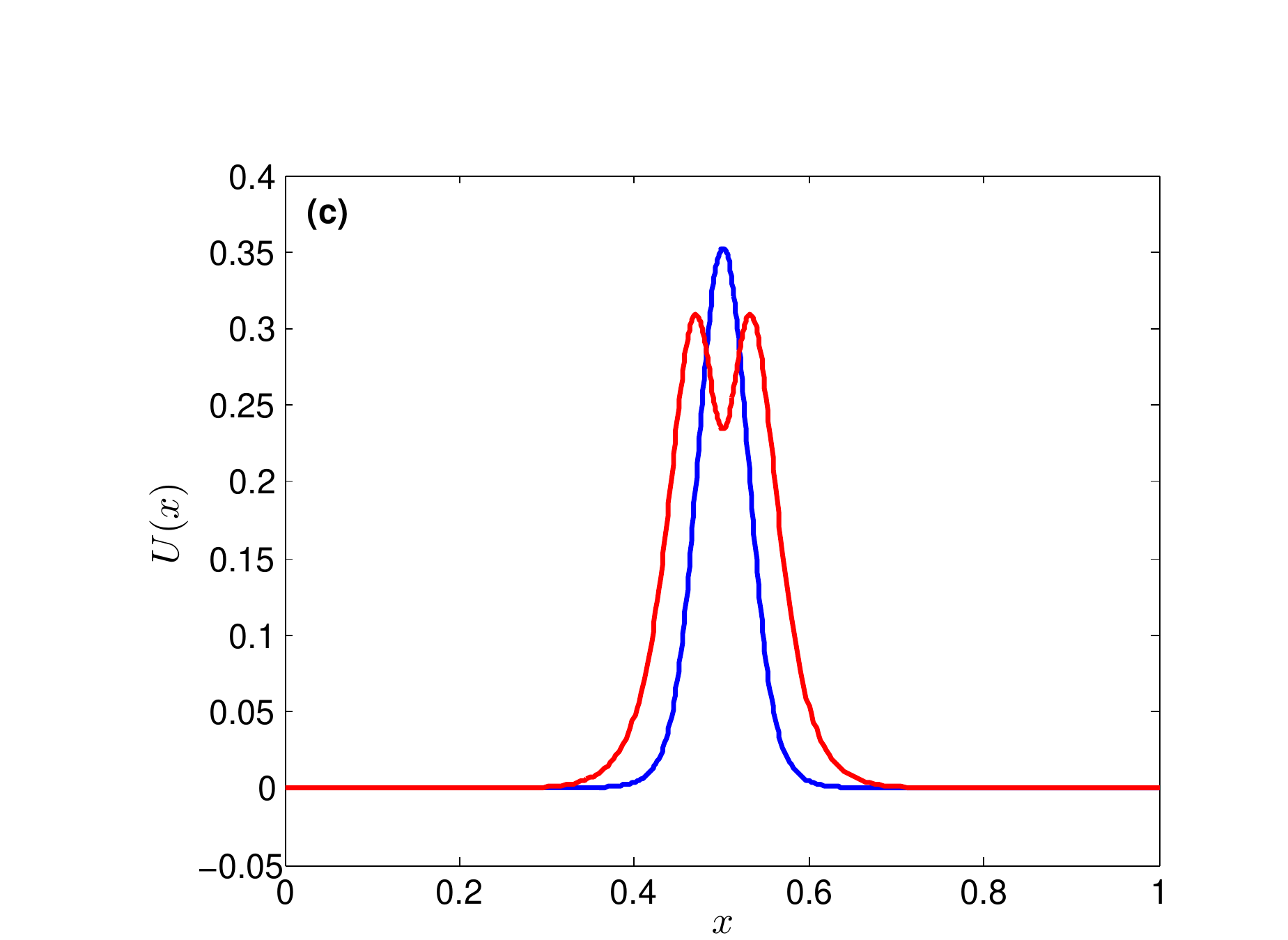}\\
\end{minipage}
&  
\begin{minipage}{.4\textwidth}
       \includegraphics[trim=.2cm 9cm 2cm 7cm, clip=true, totalheight=0.2\textheight]{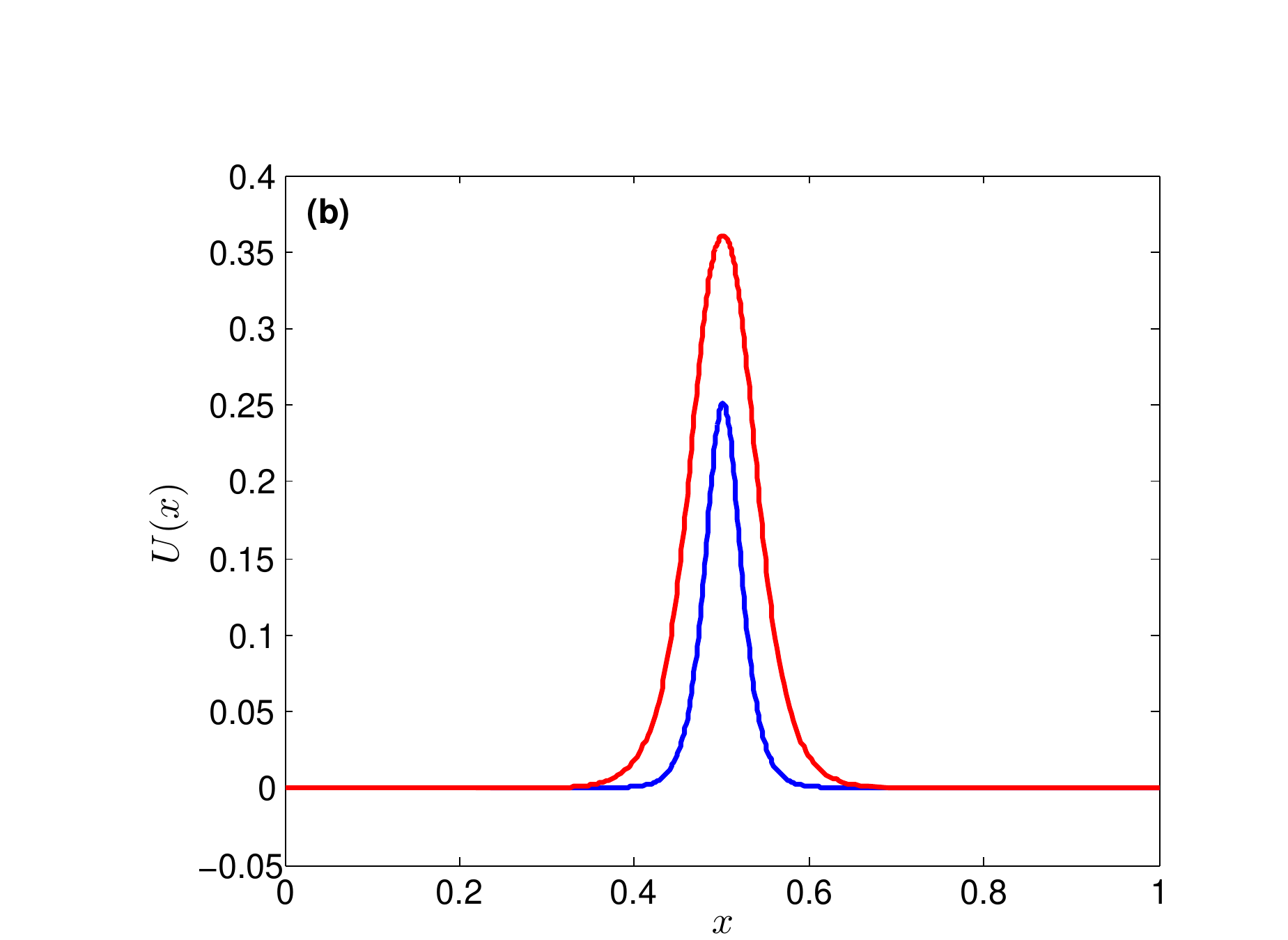}\\
         \includegraphics[trim=.2cm 9cm 2cm 7cm, clip=true, totalheight=0.2\textheight]{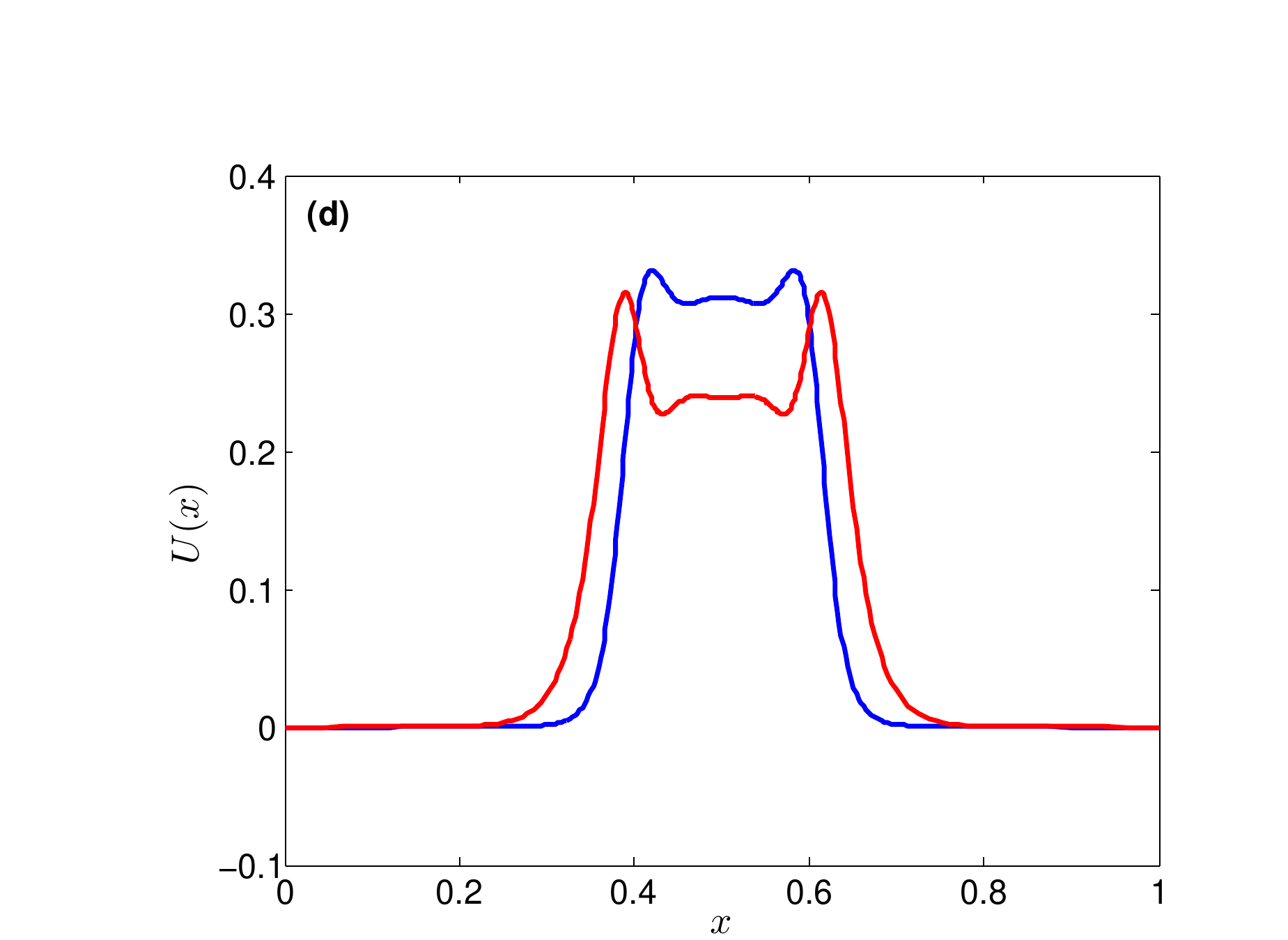}\\
\end{minipage}
  \end{tabular}
  
\caption{Examples of solutions to  (\ref{eq:rs})  in the strong damping limit with $\epsilon=0.5$, $F=2.304$, $\mu=-0.125$, $\alpha=1$, $\beta=-2$, $\nu=2$, $\omega=1+\nu \epsilon^2$, and $C=-1-2.5i$. The bistability region is between $F_0= 2.3083$ and $F_d=1.2228$. Dotted lines in (a) represent the real (blue) and imaginary (red) parts of $U_{loc}$.}
\label{fig:strong_example}
\end{figure}
\section{Conclusion} 
In the present study we examine the possible existence of spatially localized structures in the model PDE (\ref{eq:rs}) with time dependent parametric forcing. Since bistability is known to lead to the formation of localized solutions, we consider subcritical bifurcations from the zero state. 
The localized solutions we find are time dependent, unlike most previous work on this class of problems; they oscillate with half the frequency of the driving force.
In the weak damping, weak forcing limit, the solutions and bifurcations of the PDE are accurately described by its amplitude equation, the forced complex Ginzburg-Landau (FCGL) equation. Our work uses results in \cite{BYK}, where localized solutions are observed in the FCGL equation in 1D. We reduce the FCGL equation to the Allen--Cahn equation to find an asymptotically exact spatially localized solution of the PDE analytically, close to onset.\\
By continuing the numerical solution of the PDE model (\ref{eq:rs}) that we take from time-stepping as an initial condition, we found the branch of localized states. The stability of this branch was determined by time-stepping, and the region where stable localized solutions occur was found. The saddle-node bifurcations on the snaking curve arise from pinning associated with the decaying spatial oscillations on either edge of the flat state.\\
The numerical examples we give in this paper indicate how localized solutions exist in 1D, and show excellent agreement between the PDE model and the FCGL equation. The agreement remains qualitatively good even with strong damping and strong forcing. In strong damping limit, we reduce the PDE directly to the Allen--Cahn equation analytically, close to onset. By continuing the approximate solution, examples of localized oscillons are observed numerically. 

In the current work the preferred wave number is zero, so our results are directly relevant to localized pattern found in Turing systems, such as those found in \cite{TC,VE}. In contrast, in the Faraday wave experiment, the preferred wavenumber is non-zero, and so this work is not directly relevant to the oscillons that are observed there. Our interest next is to find and analyze spatially localized oscillons with non-zero wave number in the PDE model, both in 1D and in 2D. This will indicate how localized solutions might be studied in (for example) the Zhang--Vinals model \cite{ZV}, and how the weakly nonlinear calculations of \cite{SG} might be extended to the oscillons observed in the Faraday wave experiment.
\\

\textbf{Acknowledgments}. We are grateful for interesting discussions with A.D. Dean, E. Knobloch, K. McQuighan, and B. Sandstede. We also acknowledge financial support from King Abdullah Foreign Scholarship Program.

\end{document}